\documentclass[11pt]{amsart}
\usepackage{bbm}
\usepackage{amsfonts}
\usepackage{mathrsfs}
\usepackage{amsmath}
\usepackage{enumerate}

\renewcommand{\(}{\left( }
\renewcommand{\)}{\right) }

\renewcommand{\theequation}{\theequation. \arabic{equation}}
\numberwithin{equation}{section}
\newtheorem{thm}{Theorem}[section]
\newtheorem{cor}{Corollary}[section]

\newtheorem{prop}{Proposition}[section]
\newtheorem{defn}{Definition}[section]

\newtheorem{exmp}{Example}[section]

\newtheorem{special case}{Special case}[section]

\def\squarebox#1{\hbox to #1{\hfill\vbox to #1{\vfill}}}

\begin{document}
\title[Gauss summation and Ramanujan type series]
{Gauss summation and Ramanujan type series for $1/{\pi}$}
\author{Zhi-Guo Liu  }
\address{Department of Mathematics, East China Normal University,
500 Dongchuan Road, Shanghai 200241, P.R.
China}\email{zgliu@math.ecnu.edu.cn, liuzg@hotmail.com}
\date{\today}
\thanks{This work was supported  by the National Science Foundation of China
and Shanghai Natural Science Foundation (Grant No. 10ZR1409100). }
\begin{abstract}
Using some properties of the gamma function and the well-known
Gauss summation formula for the classical hypergeometric series,
we prove a four-parameter series expansion formula,
 which can produce infinitely many Ramanujan type series for $1/\pi$.
\end{abstract}
\maketitle
\section{Introduction}
\setcounter{equation}{0}
The gamma function $\Gamma(z)$ can be defined by the formula
\cite[p. 76]{Young}
\begin{equation*}
\frac{1}{\Gamma(z)}=ze^{{\gamma}z}\prod_{n=1}^\infty \(1+\frac{z}{n}\)e^{-z/n},
\end{equation*}
where $\gamma$ is the Euler constant defined as
\[
\gamma=\lim_{n\to \infty}\(1+\frac{1}{2}+\cdots+\frac{1}{n}-\log n\).
\]

$\Gamma(z)$ is meromorphic in the entire complex plane and has
simple poles at $z=0, -1, -2, \ldots.$ It is easy to verify that $\Gamma(1)=1$ and
$\Gamma(z)$ satisfies  the recurrence relation $\Gamma(z+1)=z\Gamma(z).$ It follows
that for every positive integer $n,$ we have $\Gamma(n)=(n-1)!.$ Using the recurrence
relation for the gamma function and $\Gamma(1/2)=\sqrt{\pi},$ we can find the following
proposition.
\begin{prop}\label{ppgampi} If $n$ is a  nonnegative integer, then we have
\begin{equation*}
\Gamma(n+1/2)=\frac{(2n)!}{4^n n!}\sqrt{\pi},
\ \Gamma(1/2-n)=\frac{(-1)^n 4^n n!}{(2n)!}\sqrt{\pi}.
\end{equation*}
\end{prop}

One of the most important properties of $\Gamma(z)$ is the Euler reflection formula
\cite[p. 9]{Andrews-Askey-Roy}, \cite[p. 78]{Young}.
\begin{prop}\label{ppeuler} {(\rm Euler's reflection formula).}
\[
\Gamma(z)\Gamma(1-z)=\frac{\pi}{ \sin \pi z}.
\]
\end{prop}
\begin{defn}\label{grsf}
For any complex $\alpha,$ we define the general rising shifted factorial
by
\[
 (z)_{\alpha}={\Gamma(z+\alpha)}/{\Gamma(z)}.
 \]
\end{defn}
In particular, for every non-negative integer $n, $ we have
\begin{equation}
(z)_0=1, \ (z)_n={\Gamma(z+n)}/{\Gamma(z)}=z(z+1)\cdots (z+n-1),
\label{gamEqn1}
\end{equation}
and for every positive integer $n,$
\begin{equation}
 (z)_{-n}=\frac{1}{(z-1)\cdot(z-2)\cdots(z-n)}.
 \label{agamEqn2}
\end{equation}

The well-known Gauss summation formula is stated in the following theorem
\cite[p. 66]{Andrews-Askey-Roy}, \cite[p. 102]{poole}.
\begin{thm} \label{thmgauss}
{\rm( The Gauss summation formula ).}
If $c$ is not zero or  a negative integer and $Re(c-a-b)>0,$ then we have
\[
\sum_{n=0}^\infty \frac{(a)_n (b)_n}{(c)_n n!}
=\frac{\Gamma(c)\Gamma(c-a-b)}{\Gamma(c-a)\Gamma(c-b)}.
\]
\end{thm}

In his famous paper \cite{Ram-Mod}, Ramanujan recorded a total of
$17$ series for $1/{\pi}$ without proofs.
Ramanujan's series for $1/\pi$ were not extensively studied until
around 1987. The Borwein brothers \cite{Pi-AGM} provided rigorous proofs of
Ramanujan's series  for the first time and also obtained many new series for
$1/{\pi}.$  Some remarkable extensions of them were given by the Chudnovsky
brothers \cite{Chudnovskys}.

It should be pointed out that before Ramanujan, some mathematicians had
derived some series expansions for $1/\pi$. For example, G. Bauer in 1859 \cite{Bauer}
obtained some series expansions for $1/\pi$ using Legendre polynomials,
and J. W. L. Glaisher in 1905 \cite{Glaisher} gave a systematic study
on the series expansions for $1/\pi$ using the theory of elliptic functions.

Many new Ramanujan type series for $1/{\pi}$ have been published recently,
see for example, \cite{Baruah-Berndt}, \cite{Chan-Chan-Liu},
\cite{ChuaKS}, \cite{Guillera}, \cite{Levrie}. One may consult the survey paper
\cite{Baruah-Berndt-Chan} for the interesting history of Ramanujan type
series for $1/{\pi}.$

 Very recently, Chu \cite{Chu} drived numerous Ramanujan type series for $1/{\pi}$ and $\pi$
 using Dougall's  bilateral $_2H_2$ series and the summation by parts formula.

This paper is motivated by \cite{Chu} and \cite{Levrie}.  In this paper we will use
the general rising shifted factorial and the
Gauss summation formula to prove the following four-parameter series
expansion formula for $1/\pi$.
\begin{thm}\label{rampithm} For any complex $\alpha$ and $Re(c-a-b)>0,$  we have
\begin{equation*}
\sum_{n=0}^\infty \frac{(\alpha)_{a+n}(1-\alpha)_{b+n}}{n!\Gamma(c+n+1)}
=\frac{ (\alpha)_a (1-\alpha)_b \Gamma(c-a-b)}
{ (\alpha)_{c-b}(1-\alpha)_{c-a}}
\times \frac{\sin \pi \alpha}{\pi}.
\end{equation*}
\end{thm}
 When $a, b$ and $c$ are positive integers, it is obvious that every term of the series on the left
 hand side of the above equation is a  rational function of $n.$
 Hence Theorem \ref{rampithm} allows us  to derive infinitely many series expansion formulas for
$1/\pi.$

The remainder of the paper is organized as follows. In Section $2,$
we prove Theorem \ref{rampithm} using Definition \ref{grsf} and
the Gauss summation formula.
Some special cases of Theorems \ref{rampithm} are given in Section $3.$
\section{proof of Theorem \ref{rampithm}}
\setcounter{equation}{0}
Using (\ref{gamEqn1}) we find that the identity in Theorem \ref{thmgauss}
can be rewritten as
\[
\sum_{n=0}^\infty \frac{\Gamma(a+n)\Gamma(b+n)}{n! \Gamma(c+n)}
=\frac{\Gamma(a)\Gamma(b)\Gamma(c-a-b)}{\Gamma(c-a)\Gamma(c-b)}.
\]
Replacing $(a, b, c)$ by $(a+\alpha, b+1-\alpha, c+1)$ in the above equation,
we find that
\begin{equation}
\sum_{n=0}^\infty \frac{\Gamma(a+n+\alpha)\Gamma(b+n+1-\alpha)}{n! \Gamma(c+n+1)}
=\frac{\Gamma(a+\alpha)\Gamma(b+1-\alpha)\Gamma(c-a-b)}{\Gamma(c-a+1-\alpha)\Gamma(c-b+\alpha)}.
\label{Eqn1}
\end{equation}
Using the general rising shifted factorial in Definition \ref{grsf}, we easily find that
\begin{equation*}
\Gamma(a+\alpha)=(\alpha)_{a} \Gamma(\alpha), \ \Gamma(b+1-\alpha)=(1-\alpha)_{b} \Gamma(1-\alpha),
\end{equation*}
\begin{equation*}
\Gamma(a+n+\alpha)=(\alpha)_{a+n} \Gamma(\alpha), \ \Gamma(b+n+1-\alpha)=(1-\alpha)_{b+n} \Gamma(1-\alpha),
\end{equation*}
\begin{equation*}
\Gamma(c-a+1-\alpha)=(1-\alpha)_{c-a} \Gamma(1-\alpha), \ \Gamma(c-b+\alpha)=(\alpha)_{c-b} \Gamma(\alpha).
\end{equation*}
Substituting these equations into (\ref{Eqn1}) and dividing both sides by $\Gamma(\alpha)\Gamma(1-\alpha)$,
we find that
\begin{equation*}
\sum_{n=0}^\infty \frac{(\alpha)_{a+n}(1-\alpha)_{b+n}}{n!\Gamma(c+n+1)}
=\frac{ (\alpha)_a (1-\alpha)_b \Gamma(c-a-b)}
{\Gamma(\alpha)\Gamma(1-\alpha) (\alpha)_{c-b}(1-\alpha)_{c-a}}
\end{equation*}
Replacing $\Gamma(\alpha)\Gamma(1-\alpha)$ by $\pi/{\sin \alpha \pi }$ in the right hand side
of the above equation, we complete the proof of Theorem \ref{rampithm}.
\section{Some special cases}
In this section we will give some interesting special cases of Theorem \ref{rampithm}.
\begin{cor}\label{acor}{\rm($\alpha=\frac{1}{2}$ in  Theorem \ref{rampithm})}.
For $Re(c-a-b)>0,$
\begin{equation*}
\sum_{n=0}^\infty \frac{(\frac{1}{2})_{a+n}(\frac{1}{2})_{b+n}}{n!\Gamma(c+n+1)}
=\frac{ (\frac{1}{2})_a (\frac{1}{2})_b \Gamma(c-a-b)}
{\pi (\frac{1}{2})_{c-a}(\frac{1}{2})_{c-b}}.
\end{equation*}
\end{cor}
\begin{special case} \label{sp1} \rm($a=b=0$ and $c=k$ in Corollary \ref{acor}).
If $k$ is a positive integer, then
\[
\frac{(k-1)!}{\pi (\frac{1}{2})_k^2}=\sum_{n=0}^\infty \frac{(\frac{1}{2})_n^2}{n!(k+n)!}.
\]
\end{special case}
\begin{exmp}\label{aexmp1} \rm($k=1$ in Special case \ref{sp1}: Glaisher \cite[p. 174]{Glaisher}).
\begin{equation*}
\frac{4}{\pi}=\sum_{n=0}^\infty \frac{(\frac{1}{2})_n^2}{(n+1)n!^2}.
\end{equation*}
\end{exmp}
\begin{exmp}\label{aexmp2} \rm($k=2$ in Special case \ref{sp1}).
\begin{equation*}
\frac{16}{9\pi}=\sum_{n=0}^\infty \frac{(\frac{1}{2})_n^2}{(n+1)(n+2)n!^2}.
\end{equation*}
\end{exmp}
\begin{special case}\label{sp2} \rm($a=b=-1$ and $c=k$ in Corollary \ref{acor}).
If $k$ is a nonnegative integer, then
\[
\frac{4(k+1)!^2}{\pi (\frac{1}{2})_{k+1}^2}
=4k+5+(k+1)!\sum_{n=1}^\infty \frac{(\frac{1}{2})_n^2}{(n+1)!(k+n+1)!}.
\]
\end{special case}
\begin{proof}
If $(a, b, c)=(-1, -1, k)$, then Corollary \ref{acor} becomes
\begin{align*}
\frac{ (\frac{1}{2})_{-1}^2 \Gamma(k+2)}
{\pi (\frac{1}{2})_{k+1}^2}
=\frac{(\frac{1}{2})_{-1}^2}{k!}+\frac{1}{(k+1)!}
+\sum_{n=2}^\infty \frac{(\frac{1}{2})_{n-1}^2}{n!\Gamma(k+n+1)}.
\end{align*}
Using (\ref{agamEqn2}), we easily find $({1}/{2})_{-1}=-2.$  Thus the above equation
becomes
\begin{align*}
\frac{ 4(k+1)!}
{\pi (\frac{1}{2})_{k+1}^2}
=\frac{4}{k!}+\frac{1}{(k+1)!}
+\sum_{n=2}^\infty \frac{(\frac{1}{2})_{n-1}^2}{n!\Gamma(k+n+1)}.
\end{align*}
Multiplying both sides by $(k+1)!$,  and making the variable change $n \to n+1,$
we arrive at Special case \ref{sp2}.
\end{proof}
\begin{exmp}\label{aexmp3} \rm ($k=0$ in Special case \ref{sp2}: Glaisher\cite[p. 174]{Glaisher}).
\begin{equation*}
\frac{16}{\pi}
=5+\sum_{n=1}^\infty \frac{(\frac{1}{2})_n^2}{(n+1)!^2}.
\end{equation*}
\end{exmp}
\begin{exmp}\label{aexmp4} \rm ($k=1$ in Special case \ref{sp2}).
\begin{equation*}
\frac{256}{9\pi}
=9+2\sum_{n=1}^\infty \frac{(\frac{1}{2})_n^2}{(n+2)(n+1)!^2}.
\end{equation*}
\end{exmp}

\begin{special case} \label{sp3}\rm ($a=b=-2$ and $c=k$ in Corollary \ref{acor}).
If $k$ is a nonnegative integer, then
\[
\frac{32(k+2)!(k+3)!}{\pi (\frac{1}{2})_{k+2}^2}
=32k^2+168k+217+18(k+2)!\sum_{n=1}^\infty \frac{(\frac{1}{2})_{n}^2}
{(n+2)!(k+n+2)!}.
\]
\end{special case}
\begin{proof}
If $a=b=-2$ and $c=k,$ then  Corollary \ref{acor} becomes
\begin{align*}
\frac{ (\frac{1}{2})_{-2}^2 \Gamma(k+4)}
{\pi (\frac{1}{2})_{k+2}^2}
=\frac{(\frac{1}{2})_{-2}^2}{k!}+\frac{(\frac{1}{2})_{-1}^2}{(k+1)!}
+\frac{1}{2(k+2)!}
+\sum_{n=3}^\infty \frac{(\frac{1}{2})_{n-2}^2}{n!\Gamma(k+n+1)}.
\end{align*}
Using (\ref{agamEqn2}),  we find that
$(1/2)_{-2}=4/3.$ Thus we conclude that
\begin{align*}
\frac{16(k+3)!}
{9\pi (\frac{1}{2})_{k+2}^2}
=\frac{16}{9k!}+\frac{4}{(k+1)!}
+\frac{1}{2(k+2)!}
+\sum_{n=3}^\infty \frac{(\frac{1}{2})_{n-2}^2}{n!\Gamma(k+n+1)}.
\end{align*}
Multiplying both sides of the above equation by $18(k+2)!$,
we arrive at Special case \ref{sp3}.
\end{proof}
\begin{exmp} \label{aexmp5} \rm ($k=0$ in Special case \ref{sp3}).
\begin{equation*}
\frac{2048}{3\pi}=217+36\sum_{n=1}^\infty \frac{(\frac{1}{2})_{n}^2}
{(n+2)!^2}.
\end{equation*}
\end{exmp}
\begin{cor}\label{bcor}{\rm($\alpha=\frac{1}{3}$ in  Theorem \ref{rampithm})}.
For $Re(c-a-b)>0,$
\begin{equation*}
\sum_{n=0}^\infty \frac{(\frac{1}{3})_{a+n}(\frac{2}{3})_{b+n}}{n!\Gamma(c+n+1)}
=\frac{ \sqrt{3}(\frac{1}{3})_a (\frac{2}{3})_b \Gamma(c-a-b)}
{2\pi (\frac{2}{3})_{c-a}(\frac{1}{3})_{c-b}}.
\end{equation*}
\end{cor}
\begin{special case}\label{sp4} \rm($a=b=-1$ and $c=k$ in Corollary \ref{bcor}).
If $k\ge 0$ is an integer, then
\[
\frac{\sqrt{3}k!(k+1)!}{2\pi (\frac{1}{3})_{k+1}(\frac{2}{3})_{k+1} }
=1+\frac{2}{9}k! \sum_{n=0}^\infty \frac{(\frac{1}{3})_{n}(\frac{2}{3})_{n}}{(n+1)!(n+k+1)!}.
\]
\end{special case}
\begin{proof}
If $a=b=-1$ and $c=k$, then  Corollary \ref{bcor} becomes
\[
\frac{\sqrt{3}(\frac{1}{3})_{-1}(\frac{2}{3})_{-1}\Gamma(k+2)}
{2\pi (\frac{1}{3})_{k+1}(\frac{2}{3})_{k+1}}
=\frac{(\frac{1}{3})_{-1}(\frac{2}{3})_{-1}}{k!}
+\sum_{n=1}^\infty
\frac{(\frac{1}{3})_{n-1}(\frac{2}{3})_{n-1}}{n!(n+k)!}.
\]
Using (\ref{agamEqn2}), we find that $({1}/{3})_{-1}=-3/2$ and
$({2}/{3})_{-1}=-3.$ It follows that
\[
\frac{9\sqrt{3}(k+1)!}
{4\pi (\frac{1}{3})_{k+1}(\frac{2}{3})_{k+1}}
=\frac{9}{2k!}
+\sum_{n=0}^\infty
\frac{(\frac{1}{3})_{n}(\frac{2}{3})_{n}}{(n+1)!(n+k+1)!}.
\]
Multiplying both sides of the above equation by $2k!/9,$
we arrive at Special case \ref{sp4}.
\end{proof}
\begin{exmp}  \rm($k=0$ in Special case \ref{sp4}).
\begin{equation*}
\frac{9\sqrt{3}}{4\pi }
=1+\frac{2}{9}\sum_{n=0}^\infty \frac{(\frac{1}{3})_{n}(\frac{2}{3})_{n}}{(n+1)!^2}.
\end{equation*}
\end{exmp}
\begin{cor}\label{ccor}{\rm($\alpha=\frac{1}{4}$ in  Theorem \ref{rampithm})}.
For $Re(c-a-b)>0,$
\begin{equation*}
\sum_{n=0}^\infty \frac{(\frac{1}{4})_{a+n}(\frac{3}{4})_{b+n}}{n!\Gamma(c+n+1)}
=\frac{(\frac{1}{4})_a (\frac{3}{4})_b \Gamma(c-a-b)}
{\sqrt{2} \pi (\frac{3}{4})_{c-a}(\frac{1}{4})_{c-b}}.
\end{equation*}
\end{cor}

\begin{special case}\label{sp5} \rm ($a=b=-1$ and $c=k$ in Corollary \ref{ccor}).
If $k\ge 0$ is an integer, then
\[
\frac{k!(k+1)!}{\sqrt{2}\pi (\frac{1}{4})_{k+1}(\frac{3}{4})_{k+1} }
=1+\frac{3}{16}k! \sum_{n=0}^\infty \frac{(\frac{1}{4})_{n}(\frac{3}{4})_{n}}{(n+1)!(n+k+1)!}.
\]
\end{special case}
\begin{exmp} \rm($k=0$ in Special case \ref{sp5}).
\[
\frac{8\sqrt{2}}{3\pi }
=1+\frac{3}{16}\sum_{n=0}^\infty \frac{(\frac{1}{4})_{n}(\frac{3}{4})_{n}}{(n+1)!^2}.
\]
\end{exmp}
\begin{cor}\label{dcor}{\rm($\alpha=\frac{1}{6}$ in  Theorem \ref{rampithm})}.
For $Re(c-a-b)>0,$
\begin{equation*}
\sum_{n=0}^\infty \frac{(\frac{1}{6})_{a+n}(\frac{5}{6})_{b+n}}{n!\Gamma(c+n+1)}
=\frac{(\frac{1}{6})_a (\frac{5}{6})_b \Gamma(c-a-b)}
{2\pi (\frac{5}{6})_{c-a}(\frac{1}{6})_{c-b}}.
\end{equation*}
\end{cor}

\begin{special case}\label{sp6} \rm ($a=b=-1$ and $c=k$ in Corollary \ref{dcor}).
If $k\ge 0$ is an integer, then
\[
\frac{k!(k+1)!}{2\pi (\frac{1}{6})_{k+1}(\frac{5}{6})_{k+1} }
=1+\frac{5}{36}k! \sum_{n=0}^\infty \frac{(\frac{1}{6})_{n}(\frac{5}{6})_{n}}{(n+1)! (n+k+1)!}.
\]
\end{special case}
\begin{exmp} \rm($k=0$ in Special case \ref{sp6}).
\begin{equation*}
\frac{18}{5\pi }
=1+\frac{5}{36}\sum_{n=0}^\infty \frac{(\frac{1}{6})_{n}(\frac{5}{6})_{n}}{(n+1)!^2}.
\end{equation*}
\end{exmp}
\begin{cor}\label{ecor}{\rm($\alpha=\frac{1}{10}$ in  Theorem \ref{rampithm})}.
For $Re(c-a-b)>0,$
\begin{equation*}
\sum_{n=0}^\infty \frac{(\frac{1}{10})_{a+n}(\frac{9}{10})_{b+n}}{n!\Gamma(c+n+1)}
=\frac{(\sqrt{5}-1)(\frac{1}{10})_a (\frac{9}{10})_b \Gamma(c-a-b)}
{4\pi (\frac{9}{10})_{c-a}(\frac{1}{10})_{c-b}}.
\end{equation*}
\end{cor}
\begin{special case}\label{sp7} \rm ($a=b=-1$ and $c=k$ in Corollary \ref{ecor}).
 If $k\ge 0$ is an integer, then
\[
\frac{(\sqrt{5}-1)k!(k+1)!}{4\pi (\frac{1}{10})_{k+1}(\frac{9}{10})_{k+1} }
=1+\frac{9}{100}k! \sum_{n=0}^\infty \frac{(\frac{1}{10})_{n}(\frac{9}{10})_{n}}{(n+1)! (n+k+1)!}.
\]
\end{special case}

\begin{exmp}\rm($k=0$ in Special case \ref{sp7}).
\begin{equation*}
\frac{25(\sqrt{5}-1)}{9\pi}
=1+\frac{9}{100}\sum_{n=0}^\infty \frac{(\frac{1}{10})_{n}(\frac{9}{10})_{n}}{(n+1)!^2}.
\end{equation*}
\end{exmp}

\begin{cor} \label{fcor} {\rm($\alpha=\frac{1}{5}$ in  Theorem \ref{rampithm})}.
For $Re(c-a-b)>0,$
\begin{equation*}
\sum_{n=0}^\infty \frac{(\frac{1}{5})_{a+n}(\frac{4}{5})_{b+n}}{n!\Gamma(c+n+1)}
=\frac{\(\sqrt{10-2\sqrt{5}}\) (\frac{1}{5})_a (\frac{4}{5})_b \Gamma(c-a-b)}
{4\pi (\frac{4}{5})_{c-a}(\frac{1}{5})_{c-b}}.
 \end{equation*}
\end{cor}

\begin{special case}\label{sp8} \rm ($a=b=-1$ and $c=k$ in Corollary \ref{fcor}).
If $k\ge 0$ is an integer, then
\[
\frac{(\sqrt{10-2\sqrt{5}})k!(k+1)!}{4\pi (\frac{1}{5})_{k+1}(\frac{4}{5})_{k+1} }
=1+\frac{4}{25}k! \sum_{n=0}^\infty \frac{(\frac{1}{5})_{n}(\frac{4}{5})_{n}}{(n+1)! (n+k+1)!}.
\]
\end{special case}

\begin{exmp}\rm($k=0$ in Special case \ref{sp8}).
\begin{equation*}
\frac{25(\sqrt{10-2\sqrt{5}})}{16\pi}
=1+\frac{4}{25}\sum_{n=0}^\infty \frac{(\frac{1}{5})_{n}(\frac{4}{5})_{n}}{(n+1)!^2}.
\end{equation*}
\end{exmp}
Using the same method as used in this paper,  in \cite{Liu 2011}
we use Dougall's $_5F_4$ summation for the classical hypergeometric functions
to derive a general series expansion
formula which can produce infinitely many
Ramanujan type series for $1/\pi^2$.

\paragraph {\bf Acknowledgements}
The author would like to thank the anonymous referee for many invaluable suggestions.

\end{document}